\documentclass[10pt]{article}
\usepackage{array,amsmath}
\usepackage{amssymb}
\usepackage{amsthm}
\usepackage{graphicx}
\usepackage{titlesec}
\usepackage{enumerate}
\usepackage{soul}

\titleformat{\section}{}{\bf\S\thesection.}{0.5em}{}

\def\bN{\mathbb{N}}
\def\bR{\mathbb{R}}
\def\bQ{\mathbb{Q}}

\def\cP{\mathcal{P}}
\def\cC{\mathcal{C}}
\def\s*{\textsuperscript{*}}
\begin{document}
\title{A Foundation for the Core Mathematician}
\author{David Mumford and Sy-David Friedman}
\maketitle

A large majority of research mathematicians consider the issue of  the ``Foundations of Mathematics" long settled, accepting the axiom system of Zermelo-Fraenkel with Choice (ZFC)  as adequate for everything they do. Alternatives are the various intuitionist, constructivist and finitist systems as well as the more recent homotopy-theoretic univalent foundations. But most mathematicians are happy to use ZFC.  

However, the very word foundations means a firm and fixed non-contro-versial basis, which supports all further work. In particular, foundations should support the need of research mathematicians to believe that all conjectures they work on have an answer, right or wrong (the Truth-Platonist viewpoint). On the other hand, the dramatic developments in set theory regarding large cardinal axioms  (see \cite{Ka}) and Cohen''s method of forcing new models  have resulted in a proliferation of set theories and models for the universe of sets.  What set theory you adopt has significant consequences for questions in classical descriptive set theory. There is no consensus on how ``high" you can go with infinities  nor on other possible axioms\footnote{Another example is the axiom of determinacy (AD), stating that certain games involving 2 players making alternate choices from some set {\it ad infinitum} must have winning strategies, has dramatic consequences for descriptive set theory, but this is adopted by some set-theorists and rejected by others.  There are also new forcing axioms, such as PFA (the proper forcing axiom) which have dramatic consequences.} This vigorous debate over how best to extend ZFC to accommodate the developments alluded to above continues to the present day. G\"odel believed in the Platonistic view  that set theory has a unique and well-defined universe of ``real world'' sets, a view still supported by some set-theorists (\cite{Wo}), while others feel that due to its openness and tolerance of different interpretations, set theory must accept a ``multiverse'' with a range of possible universes of set theory (\cite{SFr}, \cite{AFHT}). 

For these reasons, ZFC does not provide an adequate foundation for mathematics as a whole. But what do most mathematicians deal with in their work? The set $\bR$ of real numbers occupies a special place among all sets: it is used in the real world to measure space, time, mass and derivative quantities, using a unit and taking ratios. It is not unreasonable to say that a large majority of mathematical research concerns (i) the reals and/or (ii) the integers, together with structures derived from them both. Therein lies the ``central core'' of mathematical experience. The aim of this paper is to provide a better foundation for this ``central core'' of mathematics. 
 
In \S 1, we review how far you get without set theory. ``Second order arithmetic" or $Z_2$ turns out to be sufficient for a very large part of mathematics, including the use of trees to construct all Borel sets. But it lacks, for example, the set of all Borel sets. It feels like mental gymnastics to define topological spaces via a basis for open sets instead of a set of points.

In \S 2, we introduce a major feature of our proposal, the use of intuitions based on probability theory and relativity as expressed in Freiling's axiom
$A^\bN_{null}$  \cite{CFr} for the existence of certain sets of random numbers. He found that the infinite case of his axiom contradicted the uncountable use of the axiom of choice. We point out how unrestricted choice constructs many weird unworldly sets (and magic tricks) anyway. This paper chooses Freiling's $A^\bN_{null}$ instead of uncountable choice.

In \S 3, we spell out our axioms but, as pointed out above, we need a specific model because in models all statements are either true or false. In our approach, we take the set of real numbers as a given, modeling the physical world in the time-honored way, just as the set of natural. numbers $\bN$ is given by counting.
Thus we adopt G\"odel's Platonistic view, not for arbitrary sets, but specifically for the set $\bR$ of real numbers. 
We include Freiling and countable choice.

In \S 4, we describe the proposed model, in which all sets are given by constructions from $\bN$ and $\bR$. But, as usual since G\"{o}del's work, these constructions can use transfinite induction based on ordinals. We include a quick primer on the theory of ordinals. We propose to stop higher induction at the first ``fully $\bR$-admissible ordinal'', i.e. as soon as the resulting model satisfies ZF--  (Zermelo-Fraenkel without uncountable choice or the power set axiom). This is our main proposal for a model where every statement has a ground truth: it is reasonable to believe that all its sets can be treated as ``real world'' objects. Some of the arguments needed to prove this are technical and are given in an appendix.

in \S 5, we look at some alternatives. One is the more restricted set theory of Kripke and Platek with ``bounded quantifiers'' which may be sufficient for core mathematics. The other is category theory which does not fit cleanly into ZFC. We describe the two approaches that have been suggested, to accommodate categorical approaches.

Finally in \S 6, we look at a truly remarkable model of our axioms due to Solovay \cite{So} in which $\bR$ is countable and and all its subsets are measurable, hence Freiling's axioms hold.  This construction requires higher set theory, the existence of a model with an ``inaccessible" cardinal. However, it shows that if ZFC plus the existence of an inaccessible is consistent, then so are our axioms.  Solovay's approach also leads to an alternative ``minimalist" foundation given in this section. Other models of Freiling's axioms can also be obtained without an inaccessible using iterations of random real forcing (see the Appendix). 
 
Both authors have considered the issue of the optimal model for the foundations of mathematics \cite{M1,M2}
and for set theory \cite{SFr, AFHT}, respectively, and this led to our collaboration. The present article attempts to include enough set-theoretic background to be readable by non-specialists. We will review most definitions but good references for non-specialists who seek further clarification are the books \cite{Je, De, Ba}.

\section{\bf Do we need set theory?}\label{sec:real}
As mentioned in the Introduction, 
the natural numbers and the real numbers occupy special places among all sets.
The simplest approach to a foundation for the core mathematician would be to formulate a theory with distinct 
natural number unknowns and real unknowns. This leads to so-called {\it second order arithmetic}.

First, we need to clarify the treatment of $\bR$. It's simplest to describe a real number in the interval $[0,1]$ by its binary expansion.  There is an annoying need to identify the expansions $0. \,\vec{a}\,0\,1\,1\,1\,1\cdots= 0. \,\vec{a}\,1\,0\,0\,0\,0\cdots$ where $\vec{a}$ is any finite initial sequence, but these are only a countable set of identifications. Without this identification, the set of 1's in a binary expansion gives us an arbitrary subset of $\bN$, an element of the power set $\cP(\bN)$. The set of all sequences of 0's and 1's itself is homeomorphic to the Cantor set which we will call $\cC$. As far as set theory is concerned, there is generally no substantial difference between $\bR$, $[0,1]$, $\cC$ and the power set $\cP(\bN)$, each being reducible to any one of them.

To build a theory accommodating both $\bN$ and $\bR$, $\cP(\bN)$ is the best choice. Then a basic foundation for math can be constructed using the two sorts of variables, $a \in \bN$ and $A \in \cP(\bN)$ . This theory is called $Z_2$. Its axioms are the usual axioms of Peano arithmetic dealing with $\bN$ plus 2 axioms involving $A$'s: 
\begin{align*}
&\text{Induction: } \forall A\left(\left( 0 \in A, \forall n (n \in A \rightarrow n+1 \in A )\right) \rightarrow \forall n (n \in A)\right) \\
&\text{Comprehension: Given a formula }\phi(n), \text{ then } \exists A \forall n \left(n \in A \Leftrightarrow \phi(n)\right)
\end{align*}
Here the $\phi$ in the ``comprehension axiom" is any predicate calculus formula with quantifiers over either natural numbers $a$ or sets of natural numbers $A$. {\it The remarkable fact is that almost all core math can be developed in} $Z_2$. In particular, $\bR$, its topology and algebraic operations are easy to define in this theory. 

The great advantage of this theory is that its axioms are all clear and indisputable and it has only one obvious model: the natural numbers and {\it all} their subsets. However, there are other models where one restricts the allowed subsets of $\bN$. One even has models with only a countable set of subsets of $\bN$. Also, if one restricts the set quantifiers in comprehension in various ways, one gets constructive, finitistic and predicative ``subsystems'' of $Z_2$. There is an excellent book by Stephen Simpson \cite{Si} that studies what mathematical theorems can be proved when you restrict the quantifiers of $\phi$'s in the comprehension axiom. What is astonishing is that virtually all of standard analysis, e.g. Banach spaces and PDEs, can be developed even in such weakened systems. 

For example, one can define fully separable (also called 2nd axiom of countability) topological spaces  by enumerating a countable basis of open sets. Then points define subsets $A\subset \bN$, namely the open sets in the basis that contain them. More surprising is that Borel sets in $\bR$ (or in any Polish space) can be defined in $Z_2$.  

In all ``Intro to Analysis" books, Borel sets are defined as the smallest $\sigma$-algebra containing intervals and closed under countable unions and intersections. This involves set theory and tells you nothing about what they {\it are}. In fact, set theory is not needed if one uses {\it Borel codes} that give (not uniquely) all Borel sets. The codes are trees with a unique root, a countable set of successors to each node and such that every vertical path starting at the root is finite, ending in a leaf. These are called ``well-founded trees". The whole tree has a countable number of nodes. To get a Borel set, one further specifies a rational open interval at each leaf.  To interpret such a tree as a Borel set of reals, one works the way down from the leaves, taking unions at each odd level and intersections at each even level. For any real number $x$, one proves there is a unique ``in/out" labelling of the nodes of the tree such that, at each node, its value is the ``and' or the ``or" of the nodes just above it. 
 
Apparently, all standard theorems in mathematics involving separable spaces are provable in $Z_2$, so this theory seems to be a good candidate for ``doing math".  But what is missing in $Z_2$ is any way to talk directly about subsets of $\bR$. An individual Borel subset can be defined along with its definition as a subset of $\bR$ but there is no definition of the {\it set of all Borel subsets}. Similarly, for topological spaces, there is no way to define the {\it set of all open subsets}. One needs to embed $Z_2$ into a set theory. 

\section{\bf Freiling's axiom or Uncountable Choice}
As part of our proposed axioms, we next introduce Christopher Freiling's 1986 Axioms using {\it independent 
random real numbers}.  These axioms can be derived from intuitive ideas of independent random events bolstered 
by relativity theory. However they cannot be proven using the standard Kolmogorov reduction of probability theory 
to measure theory. Moreover his strong axiom contradicts the uncountable axiom of choice, hence the title of the 
section. However his axioms can be seen as a natural and intuitively appealing extension of probability theory 
whereas uncountable choice creates wild sets unrelated to anything real.  Freiling's axioms will be a key part 
of our proposed foundations.       

In his paper \cite{CFr} , Freiling used a ``thought experiment" of selecting random real numbers by throwing 
darts. Using binary notation for points on the unit interval, the point where the dart lands may be described by 
an infinite sequence of 0's and 1's, which may also be defined by flipping a coin an infinite number of times, 
or simply by the subset of $\bN$ given by the cases where ``heads" resulted. Of course, an infinite series of  
coin flips, let alone dart throwing, are just crude metaphors for truly random sequences (see \cite{DHM}). 
Finitists might quibble over flipping a coin infinitely many times but this is no different from countable 
choice, repeating an action $\bN$ times. Maybe taking the clicks of a geiger counter would be more truly random 
but the idea of a random choice is ultimately a mathematical idealization of randomness in the real world. 
L.E.J.Brouwer, the founder of intuitionism, came closest to defining a random sequence of 0's and 1's. He 
admitted that there were infinite sequences such that, no matter how far out you go, you cannot know the rest of 
the sequence. He called this a ``free choice sequence'' or a ``lawless sequence", see \cite{He}.  The point 
where the dart lands or the result of infinitely many coin flips are two tangible ways to connect Brouwer's 
lawless sequences with the real world.               

Henceforth, we assume the concept of a random sequence 
$$x = (b_1, b_2, b_3, \cdots), b_n \in \{0,1\}$$ 
of zeros and ones makes mathematical sense.  We can split it into two or even countably many, so called {\it 
independent} subsequences, the odd elements and the evens:
\begin{align*}&x_1 = (b_1, b_3, \cdots, b_{2n+1}, \cdots) \\
&x_2 = (b_2, b_4, \cdots, b_{2n}, \cdots)
\end{align*}
Or, using distinct primes $p_k$, one defines the countable sequence of 
independent countable sequences $x_k =\{b_{p_k}^n | n \in \bN\}$. Alternatively, the two sequences $(x_1, x_2)$ 
might be the coordinates of throwing a dart at the square dart board $[0,1]\times [0,1]$. 

Now following Freiling \cite{CFr} we imagine two dart-throws, Thrower 1 and Thrower 2, each throwing darts 
randomly at the unit interval, landing on points $x_1$ and $x_2$ respectively. We assume that the throws are 
independent in the sense that each throw is made in total ignorance of the other throw. Thus after Thrower 1 
throws his dart and later learns of the result $x_2$ of Thrower 2's throw, he cannot imagine that $x_2$ equals 
$x_1$, as the throws are being made in total ignorance of each other. Generalizing this, suppose that before the 
throws are made, both throwers are aware of some function $f:[0,1] \rightarrow \mathcal{N}([0,1])$ from $[0,1]$ 
to null subsets of $[0,1]$ and the throws are made randomly in the presence of $f$; then Thrower 1 cannot 
imagine that $x_2$ belongs to the null set $f(x_1)$ which is associated to $x_1$, again because Thrower 2 made  
his throw in total ignorance of the null set $f(x_1)$. 

The only possible exception to this is if Thrower 1 believes that somehow information about his throw got passed  
to Thrower 2 without him knowing it, before Thrower 2 makes his throw. To definitively rule out such an 
''information leak'' 
let's bring in relativity theory. Write both the location {\it and} the time of each throw as $(x_1, t_1)$ 
and $(x_2,t_2)$. It's perfectly possible neither throw is in the future cone of the other. Then {\it each} 
thrower has to wait for some interval of time after he has made his throw before any information about the other 
throw is available to him (information can be transmitted at most at the speed of light). This rules out any 
possible '' information leak''. 
This whole space-time situation is illustrated in the diagram below.        
\begin{figure*}[h]
\includegraphics[width=4in]{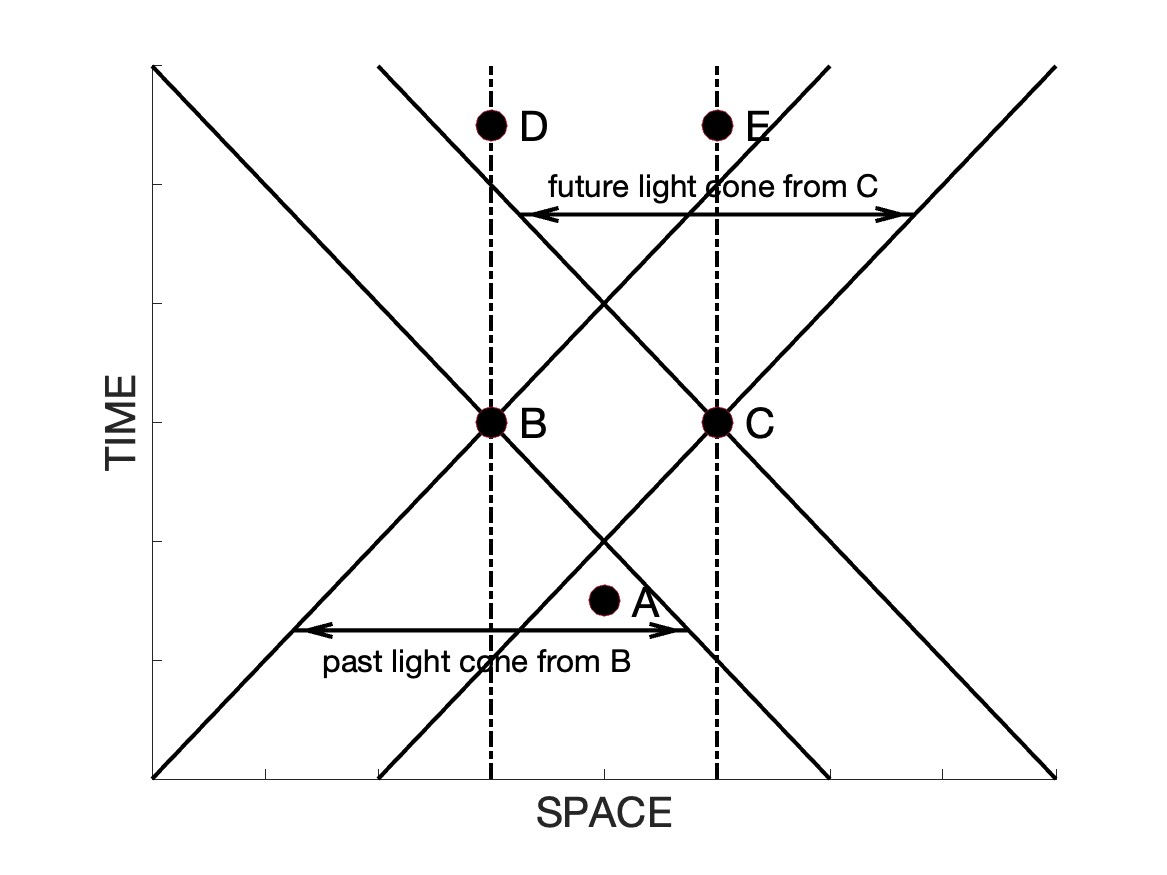}
\caption{The dotted vertical lines represent two dart throwers at different points of time, staying at two 
different fixed points of space. Each throws a dart (i.e. chooses a random real) at the space-time coordinate B and C 
respectively. The diagonal lines are light rays, the past sector (part of the full 4D past light cone) showing 
all events that can influence their throws, the future sector showing places where they can send information. 
The point A, in both past cones is where $f$ is chosen, hence at B and C, $f$ is known. B's random real is not 
$f:[0,1] \rightarrow \mathcal{N}([0,1])$,known to the dart thrower at C until, in the future, he reaches D in 
the future cone of C. Likewise the thrower 
at C learns B's throw at point E.}
\end{figure*}

So we conclude that $x_2$ is indeed not an element of $f(x_1)$, and by symmetry, neither is $x_1$ an element 
of $f(x_2)$. This is exactly Freiling's axiom:

\noindent  FREILING'S $A_{null}$; \:For all 
$f:[0,1] \rightarrow \mathcal{N}([0,1])$, 
there are (independent random) $x_1,x_2 \in [0,1]$ such that $x_1 \notin f(x_2)$ and $x_2\notin f(x_1)$ 

This situation is clarified by the key example. The continuum hypothesis
says that there is a well ordering ``$<$" of $[0,1]$ such that, for any $x_2$, the set of $x_1 < x_2$ is countable, hence a null set. This gives a wild situation: $[0,1]^2$ is divided into 3 pieces: the diagonal, $S_{12}$ where $x_1<x_2$ and its flip $S_{21}$ where $x_1 > x_2$. By symmetry, $S_{12}$ and $S_{21}$ must be the ``same" size yet every vertical line is divided into a countable part in $S_{21}$ and a set of full measure in $S_{12}$, (while every horizontal line is the opposite). This would contradict Freiling's axiom using $f(x_2)=\{x_1\big{|} x_1 < x_2\}$, i.e. $A_{null}$ \underline{disproves the continuum hypothesis.} Because G\"{o}del constructed models of ZFC where the continuum hypothesis is true, Freiling's axiom cannot be proven in ZFC.

This is a wild example. Suppose instead that the map $f$ has a nice ``graph" in the sense that the union $\Gamma_f$ of all the null sets $\{x\}\times f(x)$ is a measurable. Then Fubini's theorem applied to its indicator function shows that $\Gamma_f$ is also a null set hence random $(x_1,x_2)$ will miss both it and its flip, i.e. satisfy $A_{null}$.
 
Freiling goes further, first looking at finite sets of darts where if any one is set aside, a predicate depending on the rest is almost surely false for the one set aside. Finally, he considers an infinite sequence of dart throwers giving an independent random set of points  $\{x_1, x_2, \cdots \}$ in the interval $[0,1]$. Denote the sequence without $x_i$ by $\widehat\sigma^{(i)}$. Then:

\noindent FREILING'S  $A^\bN_{null}$: \; For all $f:[0,1] ^\bN\rightarrow \mathcal{N}$, there exists a sequence $\sigma=\{x_1, x_2, \cdots \} \in [0,1]^\bN$  such that for all $i$ , $x_i \notin f(\widehat\sigma^{(i)})$.

But the reason to accept it is that one can imagine an infinite army of dart throwers on different planets each throwing darts but at times not observable by any of the others, i.e. the times when they threw might be incomparable according to special relativity. 

Freiling then showed that, $A^\bN_{null}$ is inconsistent with the axiom of choice even if you restrict $f$ in the definition to have values in countable subsets of $[0,1]$. Hence the title of this section ``Freiling {\it or} choice". Here is Freiling's proof. Consider the tails of maps $f:\bN \rightarrow [0,1]$. By the axiom of choice, there exists a subset $Z\subset [0,1]^\bN$ with one member representing each tail. By the axiom of choice again, we can well-order $Z$ and define $f(\sigma)=$ the range of the first member of $Z$ with the same tail as $\widehat\sigma^{(1)}$.  Given $ \sigma =\{ x_1,x_2,\cdots \}$, we have that  $f(\widehat\sigma^{(1)})=f(\widehat\sigma^{(k)}) \text{ for all }k$ as $\widehat\sigma^{(1)}$ and $\widehat\sigma^{(k)}$ have a common tail. But we can choose $x_k$ in $f(\widehat\sigma^{(1)})$ as the latter contains all elements of a tail of $\widehat\sigma^{(k)}$. As  $f(\widehat\sigma^{(1)})=f(\widehat\sigma^{(k)})$, we have $x_k$ in $f(\widehat\sigma^{(k)})$, showing that $\sigma$ does not fulfill the conclusion of $A^\bN_{null}$.

His argument also shows that his axiom $A^\bN_{null}$ implies that the real numbers are not well-orderable.  Namely, if $\bR$ had a well-ordering, then so would  $\cP(\bN)$, hence the $Z$ in the above proof can be taken to consist of those sequences which do not share a tail with any sequence which is smaller 
in the sense of this well-ordering, and $Z$ inherits this well-ordering. Thus adopting the strong form of his axiom leads in a novel direction, very different from ZFC.

For those who are skeptical of dropping such a seemingly innocent axiom as the axiom of choice, there is a related consequence that may be more compelling. Thanks to Persi Diaconis for this startling result. This is a hypothetical game involving choosing infinite sequences of real numbers. You choose a sequence any way you want and put each number in a separate box. You choose an $\epsilon$ too. The magician says he will open all but one of the boxes and with probability greater than $1-\epsilon$, he will predict what is in the last box! His method can be phrased as an easy theorem that depends crucially on the same use of the axiom of choice as above, namely to fix ``tail representatives" (see \cite{Wi}, Chapter 23 for details). 

The above proof illustrates what seem to be the two common uses of the axiom of choice. One relates to a situation where a discrete group $\Gamma$ acts on a continuous object and one wants a set of representatives, one on each $\Gamma$ orbit. The other is to well-order every set so it can be measured by ordinal numbers.  Tail representatives in Freiling's argument are an instance of the first: all infinite sequences is essentially the continuous space $\bR$ and it is acted on by the discrete group of permutations of $\bN$ moving only a finite set of numbers.

Another instance is the counter-intuitive Banach-Tarski theorem that it is possible to decompose the unit ball into 5 disjoint subsets, move them around {\it by translations and rotations}, and then re-assemble them into a ball of twice the radius. The key idea, that is quite simple, is to construct a free subgroup  $\Gamma=\langle a,b \rangle \subset SO(3)$ with 2 generators. $\Gamma$ acts on the 2-sphere $S^2$ and one uses the axiom of choice to choose one point on each $\Gamma$ orbit. With a few manipulations, these define the five pieces. Yet another instance is the simplest: a set of coset representatives of $\bR/\bQ$. In all these instances, there seems to be no way to make any such choice {\it constructively}.  In the geometric cases, choice produces subsets that are unvisualizable dense clouds of points, not even fractals which typically have clear structures repeating at every scale. They have no connection to the real world or use in applied math. The main mathematical use of choice arises to deal with non-separable spaces, groups, etc. which do not seem central to the core issues of math.
 
\section{\bf Proposed Axioms}
Here are the axioms, one at a time:
 
\noindent \underline{First}, we adopt the basic uncontroversial axioms of ZFC: (i) the axiom of Extensionality -- two sets are equal if they have the same members, (ii) the axiom of foundation -- the $\in$ relation is well-founded on any set (this implies that there is no infinite chain of membership $x_{n+1} \in x_n$), (iii) the axiom of pairing -- for any two sets $x,y$ there is a set with $x$ and $y$ as its two members, (iv) axiom of union -- for any set $x$ the union of its member sets is a set, (v) the axiom of infinity -- there is a set with infinitely many members.

\noindent\underline{Second}, we want to include an axiom that says that $\bR$ exists and Freiling's axiom $A^\bN_{null}$ holds for it, based on real world understanding of randomness. 

\noindent \underline{Third}, Freiling's axiom means we must weaken the axiom of choice, and the most natural weakening is to  restrict it to countable choices.  Actually, it turns out that a stronger version is more useful. the  Axiom of Dependent Countable Choice  introduced by Paul Bernays in 1942:

\noindent {\it Axiom of Dependent Choice} {\bf (DC)}: Given a relation $R\subset X\times X$ such that $\forall a \in X, \exists b \in X \text{ with } (a,b) \in R$ and $x_0\in X$ there is a sequence $\{x_n\}$ such that $(x_n,x_{n+1}) \in R$ for all $n\geq 0$.

As a consequence, some sets may fail to have a well-ordering and the theory of cardinality is weakened: one can define $|x| \le |y|$ by the existence of an injective map from $x$ to $y$ and it does follow that $|x|\le |y| \text{ and } |y|\le |x|$ implies they are bijective. But having a surjection from $y$ to $x$ need not imply $|x| \le |y|$ and {\it two sets may have incomparable cardinalities}. There is no linearly ordered tower of cardinalities.

\noindent \underline{Fourth}, Since $\bR$ is given, the power set $\cP(\bN)$ exists but there we do not require the existence of higher power sets. In the absence of $\cP(\bR)$ and $\cP(\cP(\bN))$, we will stick with subsets of $\bR$ that are {\it constructible from $\bR$} in some sense.   Core math only uses such subsets of $\bR$ along with the higher order sets that are constructible from $\bR$ as well. Iterating the operation of taking power sets instantly produces gigantic infinities of great interest to set-theorists but without any relationship to the real world.

\noindent \underline{Fifth}, there are two other ZFC axioms that, in some form, seem central to every version of set theory:

\noindent {\it Comprehension (or Separation)}: Given a set $S$ and a formula $\phi(x)$,  then \\
$$\exists T \subset S\; \forall x \left(x \in T \Leftrightarrow \phi(x)\right)$$

\noindent{\it Collection}: Given a set $S$ and a formula $\phi(x,y)$ defining a relation $R(x,y)$ 
with domain $S$ there is a set $T$ such that for all $x\in S$ there is a $y\in T$ such that $R(x,y)$.

\noindent Taken together these two axioms imply: 

\noindent {\it Replacement}: Given a set $S$ and a formula $\phi(x,y)$ defining a function $y=f(x)$ with domain $S$, then $f(S)$ is a set.

Abbreviating ZF minus the power set axiom as ZF--, we can describe this set of axioms as T = (ZF--) + given $\bR$ +  $A^\bN_{null}$ + DC. The axiom system KP, proposed by Kripke and Platek (see \cite{De}, Ch.1, \S 11 or \cite{Ba}) has substantial overlap with the above. It was proposed as a fully predicative set theory and also rejects all choice and power sets . For this paper, the relevant part is that it restricts both comprehension and collection by requiring that the formulas $\phi$ have bounded quantifiers, i.e. for each quantifier, there is a set $w$ for which the quantifier is either $\forall x\in w$ or  $\exists x \in w$. Such formulas are referred to as $\Delta_0$ formulas.

\section{\bf Ordinals, the Proposed Model and Subsets of $\bR$}
Ordinal numbers are the backbone of set theory, a linear tower of sets on which all sets can be hung. Let's recall what ordinals are. Von Neumann invented the simplest way to define natural numbers as sets: each number is the set of all smaller numbers including zero. Thus 0 is the empty set $\emptyset$, 1 is the set whose only member is 0, 2 is the set whose members are 0 and 1, etc. $n+1=\{ 0,1,2,\cdots,n\}$. The beauty of this is that now we can name the infinite set of all natural numbers. In the context of ordinals, this is called $\omega=\{ 0,1,2,\cdots\}$, the sequence of all natural numbers with no limit. But then it goes on if you set $\omega +1 = \{ 0,1,2,\cdots, \omega\}, \: \omega+2=\{0,1,2,\cdots, \omega, \omega+1\}, \omega+\omega = \{0,1,\cdots,\omega, \omega +1, \cdots\}$ (also called $\omega\cdot 2$), etc.. This can be continued ``indefinitely", by taking the union of any $\alpha$ with the singleton set $\{\alpha\}$ to get $\alpha +1$ or, at the end of a sequence, defining a ``limit ordinal" as the set of all smaller ordinals. The sets you get in this way are all linearly ordered by $\in$ and are characterized by the additional property that any member of an ordinal is also a subset of that ordinal. Starting at any ordinal $\alpha_1$, there is no infinite descending sequence of ordinals below it $\{\alpha_1 > \alpha_2 > \alpha_3 > \cdots\}$ because $\in$ is well-founded.

Veblen (\cite{Ve}) in 1908 was the first to try to name as many of the ordinals as possible and this has been a popular game ever since. By coding tricks, many ordinals, though infinite, are still countable. For instance a bijection $f: 2\omega \xrightarrow{\simeq} \bN$ is given by $k \mapsto 2k; \;\omega+k\mapsto 2k+1$). The sup of all ordinals $k\omega$ is called $\omega^2$ and is countable by running over $\bN\times\bN$, diagonal by diagonal. Enumerating its  members, every countable ordinal ($\{\alpha_n\}|_{n\in \bN},<_\alpha$) can be described by a subset of $\bN \times \bN$, the set of all $(n,m)$ such that $\alpha_n <_\alpha \alpha_m$. The set of all countable ordinals is now the first uncountable ordinal, called $\omega_1$.  The list goes on, leading to all  the vast infinities that set-theorists have devised. 

It was G\"{o}del who realized that, parallel to the tower of ordinals, there is an increasing tower of ``constructible" sets, $L_{\alpha}$ one set for each ordinal. Here $L_{\alpha +1}$ is defined as the sets ``definable" from those in $L_\alpha$ and at limit ordinals, just taking the union of all previous constructible sets. Definable means definable by a predicate calculus formula $\phi(x)$ involving quantifiers and constants from $L_\alpha$ and one free variable $x$. The union of the $L_\alpha$'s is called simply $L$. $L$  is the simplest transitive model for every set theory $V$. Transitive models like $L$ that have the ordinals as $V$, are called ``inner models". This construction also works if some set $S$ is, for some reason ``known" (like $\bR$) and thus assumed to be given. This gives $L_\alpha(S)$ starting with $L_0(S)=S$. 

Our model is defined using $\bR$, assumed to be known, and sets constructible from $\bR$, called $\bR$-constructible sets, defined using the tower of $L_\alpha(\bR)$'s, now allowing at each stage quantifiers over $\bR$ because $\bR = L_0(\bR)$. Equivalently, we could define it using $\cP(\bN)$. It's important to realize that, in many models of ZFC, the set $\bR$ is a strict subset the standard $\bR$, e.g. it might consist only in constructible real numbers. In order to be agnostic, we take $\bR$ simply as a given, and in \S 6, we look at an alternative to the standard reals. The proposed model then is:
$$\boxed{M = L_{\alpha^*}(\bR) = L_{\alpha^*}(\cP(\bN))}$$
for a suitable ordinal $\alpha^*$. This means all sets are natural numbers, real numbers or sets constructible from $\bR$ by trans-finite induction up to $\alpha^*$. When $M$  satisfies the axioms of ZF-power, the corresponding $\alpha^*$ is called {\it fully $\bR$-admissible}. 

We need to choose $\alpha^*$ to be fully $\bR$-admissible and that $M$ satisfies DC as well. Replacement must hold in our model $L_{\alpha^*}(\bR)$, and this forces $\alpha^*$  to be fairly big. Countable ordinals $\alpha$ are all less than $\alpha^*$ because any countable ordinal can be realized (but not uniquely) as a well-ordering $<_{\alpha}$ of $\bN$:
$$\ll_\alpha=\{ (n,m) |  n <_{\alpha} m\} \subset \bN \times \bN.$$ 
By the usual correspondences, every such well-ordering can be coded as a subset of $\bN$ or as a real number in $\bR$. The conditions that a subset $\ll \subset \bN \times\bN$ is a well-ordering are (i) that it defines a linear ordering and (ii) there is no map $\bN \rightarrow \bN$ that is strictly decreasing in $\ll$ and these are expressible by simple predicates. Thus the set of well-orderings codes for $\bN$ is an $\bR$-constructible subset $WO\subset \bR$, in fact $WO \in L_1(\bR)$. But now we have the definable function $f$ that assigns to each element of $WO$ its ordinal length, so by Replacement we get the range of $f$, which is the set of all countable ordinals, i.e. $\omega_1$. More explicitly, this representation of countable ordinals is unique up to a permutation $\pi$ in the group  of all permutations $\Sigma_\bN$ of $\bN$, so 
$WO/\Sigma_\bN$
is in 1-1 correspondence with $\omega_1$. 
Thus $\omega_1$ must belong to our model $L_{\alpha^*}(\bR)$. We can't stop there: $\omega_1+1, \omega_1\cdot 2, \omega_1^2$, etc. are all definable.

But there is one curious result here. One has a surjection $\bR \twoheadrightarrow \omega_1$ obtained by sending an element of $WO$ to the ordinal length of the well-ordering it codes and other reals to $0$;  but there is no bijection because of Freiling's axiom. Thus the cardinal $\omega_1$ is in a limbo, it has no clear connection with $\bR$ in terms of cardinality\footnote{In some models of our axioms there is an injection of $\omega_1$ into $\bR$ and in other models there is none. See the Appendix.}. 
 
Our proposal for how high in the tower of $L_\alpha(R)$'s to go, can be described by looking at the sets of reals defined at each stage, $L_\alpha(\bR) \cap \cP(\bR)$. We want this to include any subset that conceivably might arise in mathematical work. We have constructed $\omega_1$ so we need to go at least to uncountable ordinals. Using the abbreviation ZF-- for Zermelo-Fraenkel without the power axiom, we denote our limit by $\alpha^*=\alpha^{\text{ZF--}}$ and define it to be the first time no new subsets of $\bR$ appear:
$$L_{(\alpha^{\text{ZF--}})}(\bR) \cap \cP(\bR) = L_{{(\alpha^{\text{ZF--}})}+1}(\bR) \cap \cP(\bR)$$

We need to prove both that there is an $\alpha$ where no new sets of reals occur and that at that point, one has a model of $T$. For the existence of any model, the axioms must be consistent, so some consistency assumption will also be needed to get a model of our theory $T$. There are definitely weaker assumptions sufficient to prove this but arguably the simplest one to state is to use the consistency of standard ZFC plus existence of an inaccessible cardinal. The consistency of this is universally accepted by set-theorists.What is ``inaccessible"? Roughly ``You cannot get there by picking yourself up by your own bootstraps". More precisely, you cannot reach it by taking power sets or by taking the supremum of a set $S$ of smaller cardinals where $S$ itself has smaller cardinality. 
Assuming this consistency, we get a model $L_{\alpha^{\text ZF-}}(\bR)$ of our theory using a  model of Solovay (\cite{So}) where all sets of reals are Lebesgue measurable.
The proof is provided in the Appendix. 

\section{\bf Quantifiers, KP and Categories}

The range of the quantifiers in predicate calculus formulas comes up not only in set theory but in mathematical practice. At one end, there is Kripke-Platek (KP) set theory described above.  In KP, both Comprehension and Collection are restricted, meaning that  the formulas involved have {\it bounded quantifiers}, called $\Delta_0$, We suspect that all formulas used in mainstream math can be made $\Delta_0$ because usually some known big set of things contains all the objects involved in each piece of work. For instance, manifolds are described by coordinate charts and gluing maps, all $\bR$-constructible sets. 

If the full force of ZF-- is not needed, one doesn't have to go quite so high in the tower of ordinals. Using definability by bounded quantifiers, the set of $\bR$-constructible subsets of $\bR$ stabilizes earlier at the ordinal $\alpha^{\text{KP}}$. This is called ``the first $\bR$-admissible ordinal (not {\it fully} $\bR$-admissible). The subsets of $\bR$ that one gets at this point are called the {\it hyperprojective} sets. 

At the other extreme, we find {\it Category Theory}. Originally, the theory arose in algebraic topology to treat spaces modulo homotopy equivalence as the objects, related by continuous maps mod homotopy. But then it became a higher level of abstraction for all sorts of mathematical objects. Thus categories $\cC$ were born: a category has two parts, first a set of objects  $S\in \text{Ob}(\cC)$ (with no internal. structure) plus a set of sets of called morphisms Hom$(S,T)$ between any two of them that can be composed satisfying the associative law. Perhaps the simplest example is the category of groups. Its objects are all groups which are sets with composition subject to the usual rules, and the usual homomorphisms between them. But this {\it cannot be a set}: like ``the set of all sets" it is too big and the comprehension axiom leads to the contradictions like ``the set of all sets that are not members of themselves".  Nonetheless ``Let $\mathcal Grp$ be the category whose objects are {\it  all groups}" is common in textbooks. This allows authors to state universal properties of all groups succinctly, e.g. the cosets of a normal subgroup form a group. 

The simplest solution is to use {\it classes}.  Recall that classes are subsets of $V$ defined by some property.  Its members are sets but classes cannot be members of each other.  In particular, the ``set of all groups" is not a set but a class. This was the approach in MacLane's 1971 book \cite{Ma}.

However, as the theory developed, it entertained ideas such as the category of all categories plus functors between them which clearly is not going to work if the categories are classes. One needs something like ``hyperclasses''. A more flexible approach is the idea of using a model $V_0 \subset V$ where $V_0$ is a set in $V$ that is called the ``universe" and has fewer ordinals than $V$. Then objects in $V_0$ are called small sets, the rest are called large sets. See \cite{MM} for details.

\section{\bf A Minimalist Foundation}
 
A key idea of contemporary set theory is the investigation of {\it countable} transitive models $M$ of set 
theory. (Transitive means the model M is a collection of sets such that, given any set in M, its members are 
also in M.) For these models, the set of  real numbers in $M$, $\bR^M = \bR \cap M$, is therefore countable, yet 
this countability is undetectable by the internal set theory of $M$. There is a fundamental theorem of Skolem 
and L\"{o}wenheim that any consistent theory in predicate calculus has a countable model. But this isn't enough 
to define such a transitive model $M$ because $M$ may have infinite descending sequences $\{\cdots \in x_3 \in 
x_2 \in x_1\}$, i.e. be ill-founded. Constructing a minimal transitive model $M$ of set theory is due to 
Shepherdson \cite{Sh} and Cohen's \cite{Co} independently. The minimal model is $L_\alpha$, where $\alpha$ is  
the least fully admissible ordinal such that $L_\alpha$ also satisfies the power set axiom, which they showed is 
countable.  The important thing to realize is that a model can be viewed both ``from the outside" and ``within 
itself".  Thus it is possible that the set of reals can be countable from the outside but of very large 
cardinality within the model.              

We seek a countable transitive model that satisfies Freiling's axiom which is in some sense ``minimal''. This 
can be done using a countable transitive model in which {\it all subsets of $\bR$ are measurable}. Although 
sounding impossible at first, such a model was discovered by Robert Solovay \cite{So}. This model satisfies all 
our axioms and certainly simplifies the behavior of sets of real numbers (e.g. every such set also has the 
property of Baire).         

Using Solovay's theory, we will construct both an extension T\s* of our previous axioms T and a ``minimal'' 
countable transitive model $N$ of T\s*.  Moreover, the theory T\s* has a strong form of ``completeness'': Any 
two {\it well-founded} models of T\s* satisfy the same sentences, and therefore incompleteness resides only in 
ill-founded models. This is an antidote to G\"odel for the mathematicians who are willing to deny the existence 
of large sets and still want there to be unique answers to all of their questions..        

Avoiding all technicalities, here is a sketch of the key idea in Solovay's theory, but in a slightly modified 
setting to get minimality. You must start by assuming ZF-- plus the existence of an inaccessible cardinal (this 
axiom is called I) is consistent. An inaccessible cardinal is just a cardinal $\aleph$ that cannot be reached by 
taking power sets and such that, for all subsets $S \subset \aleph$ where $S$ has smaller 
cardinality, their union still has smaller cardinality. Its existence is the mildest of all the infinity axioms 
proposed for set theory and is universally considered clearly consistent. Then there is a {\it countable 
model $M$} of ZF-- + I; for simplicity of exposition we assume that $M$ is also transitive. 

Now, inside this model the tower of ordinals starts as usual with finite numbers up to $\omega$, goes through 
many intermediates, past the least ordinal of size $\bR \cap M$, past the inaccessible $\kappa$ but stopping 
well before reaching another inaccessible. Internally, this is a huge tower but, viewed externally, it is 
actually countable!  Azriel L\'{e}vy first proposed using Paul Cohen's forcing to make a model that, roughly 
speaking, tells some sets, uncountable within the model, that they are actually countable. In the language of 
forcing, Solovay goes further and defines an extension $M \subset M[G]$ where $G$ is a set of surjective maps 
$\phi_\alpha: \omega \rightarrow \alpha$, one for each $\alpha < \kappa$. $M[G]$ is still countable as seen from 
the outside and has the same ordinals as $M$ but now, like an accordion opening, the position of $\omega_1$ has 
shifted way up to $\kappa$ while the rest of the tower continues up through ordinals still viewed as uncountable 
from inside. One gets new real numbers in $M[G]$, first from the sets             
$$WO_\alpha = \{\langle n,m \rangle|\phi_\alpha(n) < \phi_\alpha(m)\} \subset \omega \times \omega, \text{ for 
all }\alpha < \kappa$$      
that give explicit codes for each newly countable ordinal, second from ordering the codes: $WO^+_\alpha = 
\{WO_\beta  | \beta \leq \alpha\}$, a countable set of reals, codeable as a single real. The rest of the reals 
in $\bR^{M[G]}$ are the reals constructible relative to these. These real numbers can be divided into two parts: 
a null set equal to the union of the Borel null sets coded in M and the rest called ``random reals", not 
belonging to any Borel null set coded in $M$. In this way his random reals may be viewed as a set theory analog 
of the random reals discussed in \S 2.         
 
The least transitive model $M$ of ZF-- + I is, in fact, $L_{\alpha^{min}}$ where for no $\alpha <       
\alpha^{min}$ 
is $L_\alpha$ a model of ZF-- + I. $\alpha^{min}$ is a countable ordinal viewed from the outside, but it is 
bigger than $\kappa$, which is inaccessible inside the model $M$. Then we imitate what we did 
above in \S 4  and define the 
final minimal model $N$ via Solovay's construction as: $$N = L_{\alpha^{min}}(\bR \cap M[G],$$
the real world real numbers $\bR$ being replaced by a countable subset of them.  In this model all  
subsets of its real numbers can be proved to be measurable. But note that this model is not unique because there 
are many generic sets  $G$ corresponding to different surjections $\phi$. All these $N$'s  are the sought for 
``minimal'' models of T, one for each choice of $G$.  They can be characterized among models by two new axioms 
which say no further sets are allowed. Without going into details, here is axiom set T\s*:     
 
\begin{quote}
  1. The minimality property: there is no $\alpha$ such that $L_\alpha$ satisfies 
ZFC - Power + I.    
 
 2. There is a generic $G$ such that 
 $$V=L_{\alpha^{min}}(\bR^{M[G]})$$
 \end{quote}

It can be shown that 2 can be expressed by a first-order sentence, and therefore T\s* is a first-order theory. 
Any well-founded model of T\s* must be of the form     
$L_{\alpha^{min}}(\bR^{M[G]})$ for some generic $G$ and as any two generic $G$'s ``look the same'' it follows 
that any two well-founded models of T\s* satisfy the same first-order sentences. Thus T\s* is ``complete'' for 
well-founded models.        

\section{\bf Summary}

We have offered the core mathematician an alternative foundation based on axioms that are derived from 
mathematical experience, together with the description of two models of these axioms, as definite as the 
continuum $\bR$. The axioms contradict the axiom of choice, but validate Freiling's Axiom asserting the 
existence of independent random variables. This first foundation is based on the true continuum and believed 
to be strong enough to provide core mathematicians with everything they need to carry out their work, yet 
closely tied to basic structures built from the continuum, where core mathematics take place. Finally we 
offer a second foundation, with a ``minimalist'' ontology, which offers the 
core mathematician an  
antidote to the difficulties posed by G\"odel incompleteness. 

It is quite possible that probabilistic methods like those behind Freiling's axiom may disprove the 
measurability of 
projective sets because the $\bR$ in large cardinal models may be smaller than the standard $\bR$. This 
fits well with the perspective expressed in \cite{SFr}, that there is more evidence for the 
existence of large cardinals in inner models rather than in the full set-theoreetic universe $V$ itself.

\section*{\bf{Appendix}}
\ul{(a) If the stopping point $\alpha^\text{ZF--}$ exists, then $L_{\alpha^{\text{ZF--}}}(\bR)$ satisfies all the axioms of ZF--. }

The first step is to argue by induction on $\alpha<\alpha^\text{ZF--}$ that there is a surjection of $\bR$ onto $L_\alpha(\bR)$ which is definable over $L_\alpha(\bR$): This is obvious for $\alpha=0$ as $L_0(\bR)=\bR$ and it follows easily if $\alpha=\beta+1$ is a successor ordinal as in this case it is not hard to show that there is a surjection from $L_\beta(\bR)$ onto $L_\alpha(\bR)$ gotten by enumerating all formulas with parameters that define new sets, and this surjection is definable over $L_\alpha(\bR)$; therefore,  composing surjections provides a surjection of $\bR$ onto $L_\alpha(\bR)$ definable over $L_\alpha(\bR)$. For a limit ordinal $\alpha<\alpha^*$ we first see if there is a function definable over $L_\alpha(\bR)$ which maps $\bR$ cofinally into $\alpha$; if so we can use this function to put together surjections from $\bR$ onto smaller $L_{\bar\alpha}(\bR)$'s to obtain a surjection from $\bR$ onto all of $L_\alpha(\bR)$, definable over $L_\alpha(\bR)$. If not, then for each $n$ we can build (definably over $L_\alpha(\bR)$) an $\omega$-sequence $\alpha_0 <\alpha_1 \cdots$ of ordinals less than $\alpha$ so that $L_{\bar\alpha}(\bR)$ is a $\Sigma_n$-elementary submodel of $L_\alpha(\bR)$ where the supremum of the $\alpha_n$'s equals $\bar\alpha<\alpha$. But the latter implies that any subset of $\bR$ definable over $L_\alpha(\bR)$ is already definable over some smaller $L_{\bar\alpha}(\bR)$ and therefore is an element of $L_\alpha(\bR)$, contrary to our assumption that $\alpha$ is less than $\alpha^\text{ZF--}$. 

Next we show that $L_{\alpha^\text{ZF--}}(\bR)$ satisfies ZF--. It's enough to check Comprehension and Collection in this model. If Comprehension failed, then since by the above $L_{\alpha^\text{ZF--}}(\bR)$ contains surjections from $\bR$ onto any element of $L_{\alpha^\text{ZF--}}(\bR)$ we would get a set of reals definable over $L_{\alpha^\text{ZF--}}(\bR)$ which is not an element of $L_{\alpha^\text{ZF--}}(\bR)$, contrary to the definition of $\alpha^\text{ZF--}$.

 If Collection failed then choose a set $a$ in $L_{\alpha^\text{ZF--}}(\bR)$ and a relation $R(x,y)$ with domain $a$ which is definable over $L_{\alpha^\text{ZF--}}(\bR)$. For each $x$ in $a$ let $f(x)$ be the least ordinal $\alpha$ such that $R(x,y)$ holds for some $y$ in $L_\alpha(\bR)$. If Collection fails for $R$ then we would get a function definable over $L_{\alpha^\text{ZF--}}(\bR)$ mapping $\bR$ cofinally into $\alpha^*$, by composing $f$ with a surjection of $\bR$ onto $a$ (which exists by the first part of the proof). But then as before we can put surjections together to obtain a surjection from $\bR$ onto all of $L_{\alpha^+}(\bR)$ definable over $L_{\alpha^\text{ZF--}}(\bR)$, which by an easy diagonalization gives a subset of $\bR$ definable over $L_{\alpha^*}(\bR)$ which isi not an element of $L_{\alpha^\text{ZF--}}(\bR)$, contradicting the definition of $\alpha^*$. $\Box$ 

The ordinal $\alpha^{\text{ZF--}}$ is called the least ``fully $\bR$-admissible"  ordinal and our model $M=L_{\alpha^\text{ZF--}}(\bR)$ is called the least ``fully $\bR$-admissible set''. 

\noindent \ul{(b) Assume that there is some transitive model of ZF-- containing $\bR$ as an element which satisfies our axioms. Then  $\alpha^{\text{ZF--}}$ exists and $M(\bR)=L_{\alpha^{\text{ZF--}}}(\bR)$ satisfies our axioms; in fact it is the least transitive model of our axioms. }

To see this we work in a universe $V$ where $T^+$ holds and show that $M(\bR)$ as above exists and is a model of $T$. The fact that $M(\bR)$ exists follows from the fact that there is a transitive model of ZF - Power containing $\bR$, as if $\alpha$ is the ordinal height of such a model, $\alpha$ is fully $\bR$-admissible and therefore $\alpha^*=$ the least fully $\bR$-admissible exists. We have to check that $M(\bR)$ satisfies DC and Freiling. For DC, suppose that $R(x,y)$ is a relation definable over $M(\bR)$ such that for all $x$ in $M(\bR)$ there is a $y$ in $M(\bR)$ such that $R(x,y)$ and suppose we are given an initial $x_0$ in $M(\bR)$. Choose $\alpha_0<\alpha^*$ such that $x$ belongs to $L_{\alpha_0}(\bR)$. Inductively let $\alpha_{n+1}<\alpha^*$ be least so that for all $x$ in $L_{\alpha_n}(\bR)$ there is $y$ in $L_{\alpha_{n+1}}(\bR)$ such that $R(x,y)$. By DC (in V) we can choose $(x_0,x_1,\cdots)$ such that $x_n$ belongs to $L_{\alpha_n}(\bR)$ and $R(x_n,x_{n+1})$ for each $n$. In $L_{\alpha^*}(\bR)$ we can choose a sequence $(f_0,f_1,\cdots)$ of surjections of $\bR$ onto the $L_{\alpha_n}(\bR)$'s.  Choose (in $V$) a sequence of reals $(y_0,y_1,\dots)$ such that $f_n(y_n)=x_n$ for each $n$. Then as $L_{\alpha^*}(\bR)$ contains all reals, the sequence of $y_n$'s belongs to $L_{\alpha^*}(\bR)$ and therefore so does the sequence of $x_n$'s satisfying the conclusion of DC for $R$ in $L_{\alpha^*}(\bR)$. For Freiling, simply notice again that $L_{\alpha^*}(\bR)$ has all the reals and any function $f:\bR^N$ to null sets in $L_{\alpha^*}(\bR)$ is also such a function in $V$ so by applying Freiling for $f$ in $V$ we infer it for $f$ in $L_{\alpha^*}(\bR)$ as well. As any transitive model of ZF-- containing $\bR$ has ordinal height at least that of $M(\bR)$ it follows that $M(\bR)$ is the smallest transitive model of our theory $T$ containing $\bR$. $\Box$ 

\noindent \ul{(c) If ZFC holds and $\kappa$ is inaccessible, then there is a transitive model of our theory $T$. $T$ is a subtheory of ZF $+$ AD $+$ DC, where AD denotes the axiom of determinacy. }

First note that Freiling follows from the assertion that all sets of reals are Lebesgue measurable: Suppose $f$ maps $\bR^\bN$ to null sets. For any $k$, factor $\bR^\bN$ as 
$\bR^{\bN\setminus\{k\}} \times \bR$.
For all $\vec x$ in $\bR^{\bN\setminus\{k\}}$, $f(\vec x)$ is
null. By Fubini, the set of $\vec y$ in 
$\bR^\bN$ such that $y_k$
is in 
$f({\vec y}$ without $y_k)$
is null. The union of these null subsets of $\bR^\bN$ is null and if $\vec y$ is
not in this union then it witnesses Freiling for $f$.

Now the theory ZFC $+$ ``There is an inaccessible'' proves that there is a countable transitive model $N$ of (ZFC - Power) $+$ ``There is an inaccessible''. 
Now apply the method of Solovay (\cite{So}) to 
produces a generic extension of $N$ in which ZF $+$ DC holds and all sets of reals are  Lebesgue measurable. Then the $M(\bR)$ of that model satisfies our theory. 
As Freiling follows from the assertion that all sets of reals are Lebesgue measurable and AD proves the latter assertion, our theory is a subtheory of ZF + AD + DC. 
$\Box$ 

\noindent \ul{(d) Our theory does not imply that all projective sets of reals are Lebesgue measurable. Nor does it imply that there is no injection from $\omega_1$ into $\bR$. }

We can get an alternative model of our theory by forcing  over $L$ with a coountable support iteration of random real forcing of length $\omega_1$. The homogeneity of the forcing and the fact that it produces random reals over any given real is enough to verify that Freiling holds in the $L(\bR)$ of this model. But the set of reals of $L$ forms a non-measurable projective set in this model and as $\omega_1$ equals $\omega_1^L$, there is an injection of $\omega_1$ into $\bR$ in this model as well.  $\Box$

\end{document}